\documentclass{IEEEtran4PSCC}
\ifCLASSINFOpdf
   \usepackage[pdftex]{graphicx}
\else
   \usepackage[dvips]{graphicx}
\fi
\usepackage[cmex10]{amsmath}
\usepackage{amsthm,amssymb}
\usepackage{algorithm}
\usepackage{algorithmic}

\usepackage{url}

\makeatletter
\let\old@ps@headings\ps@headings
\let\old@ps@IEEEtitlepagestyle\ps@IEEEtitlepagestyle
\def\psccfooter#1{%
    \def\ps@headings{%
        \old@ps@headings%
        \def\@oddfoot{\strut\hfill#1\hfill\strut}%
        \def\@evenfoot{\strut\hfill#1\hfill\strut}%
    }%
    \def\ps@IEEEtitlepagestyle{%
        \old@ps@IEEEtitlepagestyle%
        \def\@oddfoot{\strut\hfill#1\hfill\strut}%
        \def\@evenfoot{\strut\hfill#1\hfill\strut}%
    }%
    \ps@headings%
}
\makeatother

\theoremstyle{definition}

\theoremstyle{proposition}

\usepackage{multirow}
\usepackage[colorlinks]{hyperref}

\begin{document}
%
\title{Accelerated Computation and Tracking of\\AC Optimal Power Flow Solutions Using GPUs}

\author{
\IEEEauthorblockN{Youngdae Kim\\Kibaek Kim}
\IEEEauthorblockA{Mathematics and Computer Science Division\\
Argonne National Laboratory\\
Lemont, IL, USA\\
\{youngdae,kimk\}@anl.gov}
}


\maketitle

\begin{abstract}
We present a scalable solution method based on an alternating direction method of multipliers and graphics processing units (GPUs) for rapidly computing and tracking a solution of alternating current optimal power flow (ACOPF) problem.
Such a fast computation is particularly useful for mitigating the negative impact of frequent load and generation fluctuations on the optimal operation of a large electrical grid.
To this end, we decompose a given ACOPF problem by grid components,  resulting in a large number of small independent nonlinear nonconvex optimization subproblems.
The computation time of these subproblems is significantly accelerated by employing the massive parallel computing capability of GPUs.
In addition, the warm-start ability of our method leads  to  faster convergence, making the method particularly suitable for fast tracking of optimal solutions.
We demonstrate the performance of our method on grids having up to 70,000 buses by solving associated optimal power flow problems with both cold start and warm start.
\end{abstract}

\begin{IEEEkeywords}
alternating current optimal power flow, alternating direction method of multipliers, graphics processing unit
\end{IEEEkeywords}

\thanksto{\noindent This research was supported by the Exascale Computing Project (17-SC-20-SC),
a collaborative effort of the U.S. Department of Energy Office of Science and the National Nuclear Security Administration.
This material was based upon work supported by the U.S. Department  of Energy, Office of Science, under Contract No. DE-AC02-06CH11357.}

\section{Introduction}
\label{sec:intro}

The growing penetration of renewable and distributed energy resources  into electrical grids brings attention to distributed solution methods for solving alternating current optimal power flow (ACOPF) problems to efficiently control and optimize a potentially huge-scale grid~\cite{molzahn2017survey}.
In contrast to centralized optimization methods such as Ipopt~\cite{Wachter06}, distributed solution methods provide a scalable computation capability, making them particularly suitable for an environment where the number of control units is ever increasing.
Because computations can be performed in a distributed fashion, they are also less vulnerable to privacy threats and cyber attacks.

Another important aspect to consider for  reliable operation of the grid under such circumstances is the ability to rapidly track an optimal solution under load and generation fluctuations.
Power generation from renewable resources is intrinsically uncertain, causing potentially larger fluctuations in generation over time than those caused by traditional fuel-based generators.
Without a timely and optimal adjustment of set points of generators, these fluctuations could incur inefficiency and instability in grid operations, often represented as repeated large deviations of the system frequency from its nominal value.
In this case, quickly tracking an optimal solution from the previous solution is of central importance.
From the optimization perspective, this tracking ability is closely related to the warm-start capability of the underlying optimization methods.

Among many different distributed solution methods, the alternating direction method of multipliers (ADMM) is of particular interest because of its superior computational performance and privacy-preserving capability, as demonstrated in~\cite{Mhanna19,SunSun21,RyuKim21}.
The algorithm is also suitable for exploiting warm start for accelerated convergence.
However, the existing ADMM-based distributed solution methods in the literature work only on CPUs, without fully realizing their algorithmic potential for the massive parallel computing of the ADMM subproblems.
The main reason is  hardware limitations; for example, a quad-core CPU can run up to 8 threads, whereas a GPU can run tens of thousands of threads and more simultaneously.

Although GPUs have shown great success in accelerating machine learning algorithms such as deep neural networks, they have made little progress in advancing the algorithms for nonlinear optimization problems such as ACOPF.
One of the main reasons is the lack of an efficient linear system solver for a large-scale sparse symmetric indefinite system of equations, for example the Karush--Kuhn--Tucker system matrix.
Most nonlinear optimization algorithms require factorization and triangular solves of such a system of equations in order to compute a Newton direction.
These computations take a significant amount (more than 80\%) of the total computation time.
However, the sparsity of the system of equations and inherently sequential nature of the algorithms leave little room for accelerating their computation on GPUs, as shown in~\cite{Tasseff19,Swirydowicz21}.

In this paper we show how we may exploit GPUs to significantly accelerate the computation time of an ADMM-based distributed solution method, while guaranteeing convergence of its iterates.
Our approach is based on a combination of the two existing ADMM algorithms for ACOPF~\cite{Mhanna19,SunSun21}.
We combine the computational and theoretical merits of those approaches, respectively, in a way that we can utilize the enormous parallel computing capability of GPUs along with benefiting from good theoretical properties.

In particular we exploit the computational advantages of~\cite{Mhanna19} by decomposing an ACOPF problem into many small  subproblems that can be solved in parallel.
Many of these subproblems have a closed-form solution and thus can be easily parallelized on GPUs.
However,  other subproblems are small nonlinear nonconvex optimization problems.
To rapidly solve such problems in batch on GPUs,
we employ our recently developed batch solver, ExaTron~\cite{ExaTron}, which is capable of significantly accelerating the computation of solutions of tens of thousands and more of small nonlinear nonconvex problems on GPUs.
We introduce our reformulation and augmented Lagrangian techniques to fully utilize ExaTron for this purpose.

While~\cite{Mhanna19} does not guarantee the convergence of its ADMM iterates, we modify its computational procedure to apply recent theoretical developments of~\cite{SunSun21} to establish convergence guarantees.
The resulting ADMM algorithm has both the computational and theoretical benefits of~\cite{Mhanna19,SunSun21}, respectively, that can be efficiently implemented and executed on GPUs.
Experimental results in Section~\ref{sec:exp} demonstrate that our ADMM algorithm on GPUs is as competitive as Ipopt for solving a large-scale ACOPF problem having up to 70,000 buses.
In addition, our method takes advantage of warm starting, which  further provides accelerated convergence, making our method particularly suitable for fast tracking of optimal solutions under load and generation fluctuations.

The rest of the paper is organized as follows.
In Section~\ref{sec:formulation} we introduce our ADMM decomposition scheme including its formulations and convergence proof.
Section~\ref{sec:implementation} briefly describes the implementation of our ADMM algorithm on GPUs.
In Section~\ref{sec:exp} we demonstrate the computational performance of our method on GPUs over large-scale grids having up to 70,000 buses with both cold start and warm start.
We compare its computation time and solution quality with those obtained from Ipopt.
We conclude  in Section~\ref{sec:conclusion} with a summary and a brief look at future work for improving ADMM.
Our algorithm has been implemented in Julia~\cite{Julia-2017} using CUDA.jl~\cite{besard2018juliagpu} and is available at \url{https://github.com/exanauts/ExaAdmm.jl.git}.


\section{ADMM formulations for ACOPFs}
\label{sec:formulation}

Our ADMM formulation is a combination of the component-based decomposition~\cite{Mhanna19} and the two-level algorithm~\cite{SunSun21}.
The component-based decomposition enables us to exploit the massive parallel computing capability of GPUs by generating many small nonlinear nonconvex problems in a form that can be efficiently solved on GPUs in parallel, while the two-level algorithm provides theoretical grounds for convergence guarantees of our ADMM iterates.
The details of those two formulations can be found in~\cite{Mhanna19,SunSun21}; here we introduce their basic ideas briefly and describe the main difference between our formulation and theirs.
At the end of this section we provide a sketch of the convergence proof of our ADMM iterates.


\subsection{ACOPF formulation}
\label{subsec:acopf}

A rectangular formulation of an ACOPF problem is given in~\eqref{eq:acopf}.
ACOPF computes economically optimal set  points of generators and voltage values within their lower and upper limits that satisfy physical laws such as Ohm's law and Kirchhoff's laws.
These physical laws are represented in power flow equations~\eqref{eq:pg}--\eqref{eq:qg} and~\eqref{eq:pij}--\eqref{eq:qji}.
Because of the nonlinearity of power flow equations, ACOPF problems are nonlinear nonconvex optimization problems and thus NP-hard.
Therefore, we typically aim to compute a stationary point as a solution, which is a point satisfying the well-known first-order optimality conditions.

\begin{subequations}
\begin{align}
    &\underset{p_{g_i},q_{g_i},w_i,\theta_i,w^R_{ij},w^I_{ij}}{\text{minimize}} \; \sum_{i \in \mathcal{B}}\sum_{g_i \in \mathcal{G}_i} f_{g_i}(p_{g_i})\\
    &\text{subject to} \notag\\
    &\sum_{g_i \in \mathcal{G}_i} p_{g_i} - p_{d_i} = g_i^Sw_i + \sum_{j \in \mathcal{B}_i}p_{ij}, \; \forall i \in \mathcal{B} \label{eq:pg}\\
    &\sum_{g_i \in \mathcal{G}_i} q_{g_i} - q_{d_i} = -b_i^Sw_i + \sum_{j \in \mathcal{B}_i}q_{ij}, \; \forall i \in \mathcal{B} \label{eq:qg}\\
    &\sqrt{p_{ij}^2 + q_{ij}^2} \le \bar{r}_{ij}, \; \forall (i,j) \in \mathcal{L} \label{eq:ij-limit}\\
    &\sqrt{p_{ji}^2 + q_{ji}^2} \le \bar{r}_{ji}, \; \forall (i,j) \in \mathcal{L} \label{eq:ji-limit}\\
    & \underline{p}_{g_i} \le p_{g_i} \le \overline{p}_{g_i}, \; \forall g_i \in \mathcal{G}_i, \forall i \in \mathcal{B}\\
    & \underline{q}_{g_i} \le q_{g_i} \le \overline{q}_{g_i}, \; \forall g_i \in \mathcal{G}_i, \forall i \in \mathcal{B}\\
    & -2\pi \le \theta_i \le 2\pi, \; \forall i \in \mathcal{B}\\
    & p_{ij} = g_{ii}w_i + g_{ij}w^R_{ij} + b_{ij}w^I_{ij}, \; \forall (i,j) \in \mathcal{L} \label{eq:pij}\\
    & q_{ij} = -b_{ii}w_i - b_{ij}w^R_{ij} + g_{ij}w^I_{ij}, \; \forall (i,j) \in \mathcal{L} \label{eq:qij}\\
    & p_{ji} = g_{jj}w_j + g_{ji}w^R_{ij} - b_{ji}w^I_{ij}, \; \forall (i,j) \in \mathcal{L} \label{eq:pji}\\
    & q_{ji} = -b_{jj}w_j - b_{ji}w^R_{ij} - g_{ji}w^I_{ij}, \; \forall (i,j) \in \mathcal{L} \label{eq:qji}\\
    & (w^R_{ij})^2 + (w^I_{ij})^2 = w_iw_j, \; \forall (i,j) \in \mathcal{L} \label{eq:wij-magnitude}\\
    & \theta_i - \theta_j = \arctan(w^I_{ij}/w^R_{ij}), \; \forall (i,j) \in \mathcal{L} \label{eq:wij-angle}
\end{align}
\label{eq:acopf}
\end{subequations}

In~\eqref{eq:acopf}, we have $\theta_{ij}=(\theta_i-\theta_j)$, $(\tilde{y}_{ij}+\frac{1}{2}\tilde{y}_{ij}^{\text{Sh}})/|a_{ij}|^2=g_{ii}+jb_{ii}$, $(-\tilde{y}_{ij})/a_{ij}^*=g_{ij}+jb_{ij}$,  $(-\tilde{y}_{ij})/a_{ij}=g_{ji}+jb_{ji}$, $\left(\tilde{y}_{ij}+\frac{1}{2}\tilde{y}_{ij}^{\text{Sh}}\right)=g_{jj}+jb_{jj}$,
where $\tilde{y}_{ij}, \tilde{y}_{ij}^\text{Sh},$  and $a_{ij}$ denote a branch series admittance, a branch shunt admittance, and a turns ratio between bus $i$ (from) and bus $j$ (to), respectively;
$g_i^S$ and $b_i^S$ represent a shunt conductance and susceptance at bus $i$;
$p_{g_i}$ and $q_{g_i}$ are variables for real and reactive powers of generator $g_i$ at bus $i$;  $w^R_{ij}$ and $w^I_{ij}$ are defined to be $v_iv_j\cos\theta_{ij}$ and $v_iv_j\sin\theta_{ij}$, respectively, with $v_i$ being the voltage magnitude at bus $i$;
$w_i$ and $\theta_i$  are for squared voltage magnitude $(=v_i^2)$ and angle at bus $i$;  and
$p_{d_i}$ and $q_{d_i}$ denote real and reactive loads at bus $i$.
The notations $\mathcal{B},\mathcal{B}_i,\mathcal{G}_i$, and $\mathcal{L}$ denote
a set of bus indices, a set of bus indices connected to bus $i$, a set of generator indices at bus $i$, and a set of line indices connecting between buses, respectively.
We note that for $(i,j) \in \mathcal{L}$ $i$ is a from-bus and $j$ a to-bus.

\subsection{Basic idea of existing ADMM formulations in~\cite{Mhanna19,SunSun21}}
\label{formulation:existing-admm}

ADMM enables the decomposition of a given problem with linear coupling constraints by forming an augmented Lagrangian and performing iterates in an alternating fashion.
For example, in order to solve the following problem
\begin{equation}
    \begin{aligned}
    & \underset{x,\bar{x}}{\text{minimize}} && f(x)+g(\bar{x})\\
    & \text{subject to} && Ax+B\bar{x}=c,
    \end{aligned}
    \label{eq:admm}
\end{equation}
an augmented Lagrangian is formed as $L_\rho(x,\bar{x},y)=f(x)+g(\bar{x})+y^T(Ax+B\bar{x}-c)+(\rho/2)\|Ax+B\bar{x}-c\|_2^2$, and variables are then updated sequentially in the order of $x, \bar{x}$, and $y$. 
These alternating updates enable us to decouple the problem into subproblems involving $x$ or $\bar{x}$ variables only, facilitating parallel computing.
The decoupling of $x$ variable from $\bar{x}$ variable is a key difference from the existing augmented Lagrangian method where we perform a joint optimization over $x$ and $\bar{x}$. 

In~\cite{Mhanna19}, a given ACOPF problem~\eqref{eq:acopf} is decomposed into components subproblems, such as generators, branches, and buses.
Its basic idea stems from the observation that generators and buses are coupled through $p_{g_i}$ and $q_{g_i}$ variables and branches and buses through $p_{ij},q_{ij},p_{ji},q_{ji},w_i,\theta_i,w_j,\theta_j$ variables for a given branch $(i,j)$.
By duplicating these variables and enforcing a consensus through coupling constraints, we can reformulate the problem into an equivalent form where it can be decomposed into component subproblems.
For example, we create $p_{g_i(i)}$ and $q_{g_i(i)}$ variables, duplicates of $p_{g_i}$ and $q_{g_i}$, respectively, and add coupling constraints $p_{g_i}-p_{g_i(i)}=0$ and $q_{g_i}-q_{g_i(i)}=0$.
Similarly, we duplicate variables between branches and buses, resulting in new $p_{ij(i)}, q_{ij(i)}, p_{ji(i)}, q_{ji(i)}, w_{i(ij)}, \theta_{i(ij)}, w_{j(ij)}, \theta_{j(ij)}$ variables and coupling constraints $p_{ij}-p_{ij(i)}=0, q_{ij}-q_{ij(i)}=0, p_{ji}-p_{ji(i)}=0, q_{ji}-q_{ji(i)}=0, w_{i(ij)}-w_i=0, \theta_{i(ij)}-\theta_i=0, w_{j(ij)}-w_j=0, \theta_{j(ij)}-\theta_j=0$.

By properly assigning these original and duplicate variables into components, we can decompose the problem into component subproblems.
The generator component $g_i$ is assigned with $p_{g_i}$ and $q_{g_i}$, the bus component $i$ with $w_i, \theta_i$ and duplicate variables having subscript $(i)$, and the branch component $(i,j)$ with $p_{ij}, q_{ij}, p_{ji}, q_{ji}, w^R_{ij},w^I_{ij}$ and duplicate variables having subscript $(ij)$.
We then constitute~\eqref{eq:pg}--\eqref{eq:qg} with bus $i$ variables only and ~\eqref{eq:ij-limit}--\eqref{eq:ji-limit} and~\eqref{eq:pij}--\eqref{eq:wij-angle} with branch $(i,j)$ variables only.
If we apply ADMM to the resulting problem by identifying $x$ with generator and branch variables and $\bar{x}$ with bus variables in~\eqref{eq:admm}, one can easily verify that an ADMM iteration leads to subproblems consisting of individual components only, which can be solved in a massively parallel way (e.g.,~\cite{Mhanna19,ExaTron}).
In Section~\ref{sec:implementation} we discuss parallel solves of these types of subproblems on GPUs.

While~\cite{Mhanna19} may not guarantee the convergence of its ADMM iterates, the recent theoretical development of a two-level algorithm~\cite{SunSun21} guarantees convergence of the ADMM for ACOPF problems to a stationary point.
In terms of the problem formulation~\eqref{eq:admm}, the basic idea of the two-level scheme is to introduce an artificial variable $z$ to the coupling constraint and enforce it to have a zero value as follows.
\begin{equation}
    \begin{aligned}
    &\underset{x,\bar{x},z}{\text{minimize}} && f(x) + g(\bar{x})\\
    &\text{subject to} && Ax+B\bar{x} - c + z = 0, \quad z = 0
    \end{aligned}
    \label{eq:two-level}
\end{equation}
Clearly, \eqref{eq:admm} and~\eqref{eq:two-level} are equivalent.
When we solve~\eqref{eq:two-level}, however, an augmented Lagrangian method is applied on the $z=0$ constraint only.
In this case each iteration (called an outer iteration) of the augmented Lagrangian method consists of solving the problem $L_\beta(x,\bar{x},z)=f(x)+g(\bar{x})+\lambda^Tz + (\beta/2)\|z\|^2_2$ subject to $Ax+B\bar{x}-c+z=0$ and updating the multiplier $\beta$.
When solving $L_\beta$ under the coupling constraint, we apply an ADMM algorithm to it.
Iterations of this inner ADMM algorithm are called inner iterations, hence the name two-level algorithm.
Introducing such an artificial variable $z$ enables us to satisfy assumptions on the last block (in this case the $z$ variable) of the coupling constraints for the ADMM to converge:
it has the entire space as its domain and can always find a feasible point satisfying the coupling constraints for any given pair of $(x,\bar{x})$.
We refer to~\cite{SunSun21} for more details.

\subsection{Our ADMM formulation}
\label{subsec:our-admm}

Our ADMM formulation for solving ACOPFs is a combination of the aforementioned two ADMM formulations~\cite{Mhanna19,SunSun21} by taking their computational and theoretical merits, respectively.  In
~\cite{Mhanna19}  many independent subproblems of small sizes are provided that are amenable to exploiting the massive parallel computing capability of GPUs, while the  work in~\cite{SunSun21} lays theoretical grounds for convergence guarantees.
In Section~\ref{subsec:design-principles} we detail the reason why we may want to use the component-based decomposition scheme of~\cite{Mhanna19} for efficient computation on GPUs and borrow theories only from~\cite{SunSun21}.

To apply the theoretical framework of~\cite{SunSun21} to~\cite{Mhanna19}, we make the following changes.
First, an artificial variable is introduced for each coupling constraint to reformulate~\eqref{eq:admm} into~\eqref{eq:two-level}.
For example, $p_{g_i}-p_{g_i(i)}=0$ becomes $p_{g_i}-p_{g_i(i)}+z_{p_{g_i}}=0$ with constraint $z_{p_{g_i}}=0$.
We then apply an augmented Lagrangian on constraint $z=0$ only, where $z$ denotes a vector of the added artificial variables.
The reformulation is simple and straightforward, but it will provide us with theoretical grounds to prove the convergence of our ADMM iterates.


Second, each branch problem is reformulated in terms of voltage and slack variables, where slacks are for the line limit constraints~\eqref{eq:ij-limit}--\eqref{eq:ji-limit} in this case.
An augmented Lagrangian is then applied to the resulting formulation on the line limit constraints.
The main reason for this reformulation is to change the branch problems into a batch of bound-constrained optimization problems in order to use our GPU solver ExaTron~\cite{ExaTron} to accelerate their computation time.

Specifically, $w_i, w_j, w^R_{ij}, w^I_{ij}$ variables are replaced with expressions in terms of $v_i,v_j,\theta_i,\theta_j$ variables, where $w_i=v_i^2, w_j=v_j^2, w^R_{ij}=v_iv_j\cos\theta_{ij}$, and $w^I_{ij}=v_iv_j\sin\theta_{ij}$.
The power flow variables $p_{ij},q_{ij},p_{ji},q_{ji}$ are free variables and explicitly defined in~\eqref{eq:pij}--\eqref{eq:qji}.
In this case we can replace each occurrence of those variables in the problem with the expression on the right-hand side of~\eqref{eq:pij}--\eqref{eq:qji}.
After the replacement,~\eqref{eq:pij}--\eqref{eq:qji} as well as~\eqref{eq:wij-magnitude}--\eqref{eq:wij-angle} are no longer needed, and we therefore remove all those equality constraints from consideration.
To deal with line limits~\eqref{eq:ij-limit}--\eqref{eq:ji-limit}, we formulate an augmented Lagrangian on those constraints.
Two new slack variables $s_{ij}$ and $s_{ji}$ are introduced to transform the inequality constraints into equalities, resulting in $p_{ij}^2 + q_{ij}^2 + s_{ij}=0$ and $p_{ji}^2 + q_{ji}^2 + s_{ji}=0$ with $-\bar{r}_{ij} \le s_{ij} \le 0$ and $-\bar{r}_{ji} \le s_{ji} \le 0$.
These equality constraints are dualized, and augmented terms are added to form an augmented Lagrangian.

The resulting formulation for the ADMM iteration of a branch $(i,j)$ subproblem is a bound-constrained nonlinear nonconvex optimization problem as described in~\eqref{eq:branch}.
It involves fewer variables and constraints than does the original branch problem.
In~\eqref{eq:branch}, $\lambda$ and $\rho$ are a multiplier and a penalty parameter for the ADMM algorithm, and $\tilde{\lambda}$ and $\tilde{\rho}$ are for the application of an augmented Lagrangian algorithm for the problem.
\begin{equation}
    \begin{aligned}
    &\underset{\substack{v_i,v_j,\theta_{i(ij)},\\\theta_{j(ij)},s_{ij},s_{ji}}}{\text{minimize}} \;\; \lambda_{p_{ij}}(p_{ij}-p_{ij(i)}+z_{p_{ij}}) \\
    &\; + \frac{\rho_{p_{ij}}}{2}(p_{ij}-p_{ij(i)}+z_{p_{ij}})^2\\
    &\; +\lambda_{q_{ij}}(q_{ij}-q_{ij(i)}+z_{q_{ij}}) + \frac{\rho_{q_{ij}}}{2}(q_{ij}-q_{ij(i)}+z_{q_{ij}})^2\\
    &\; +\lambda_{p_{ji}}(p_{ji}-p_{ji(i)}+z_{p_{ji}}) + \frac{\rho_{p_{ji}}}{2}(p_{ji}-p_{ji(i)}+z_{p_{ji}})^2\\
    &\; +\lambda_{q_{ji}}(q_{ji}-q_{ji(i)}+z_{q_{ji}}) + \frac{\rho_{q_{ji}}}{2}(q_{ji}-q_{ji(i)}+z_{q_{ji}})^2\\
    &\; +\lambda_{w_i}(v_i^2 - w_i + z_{w_{i(ij)}}) + \frac{\rho_{w_i}}{2}(v_i^2 - w_i + z_{w_{i(ij)}})^2\\
    &\; +\lambda_{\theta_i}(\theta_{i(ij)} - \theta_i + z_{\theta_{i(ij)}}) + \frac{\rho_{\theta_i}}{2}(\theta_{i(ij)} - \theta_i + z_{\theta_{i(ij)}})^2\\
    &\; +\lambda_{w_j}(v_j^2 - w_j + z_{w_{j(ij)}}) + \frac{\rho_{w_j}}{2}(v_j^2 - w_j + z_{w_{j(ij)}})^2\\
    &\; +\lambda_{\theta_j}(\theta_{j(ij)} - \theta_j + z_{\theta_{j(ij)}}) + \frac{\rho_{\theta_j}}{2}(\theta_{j(ij)} - \theta_j + z_{\theta_{j(ij)}})^2\\
    &\; +\tilde{\lambda}_{s_{ij}}(p_{ij}^2+q_{ij}^2+s_{ij}) + \frac{\tilde{\rho}_{s_{ij}}}{2}(p_{ij}^2+q_{ij}^2+s_{ij})^2\\
    &\; +\tilde{\lambda}_{s_{ji}}(p_{ji}^2+q_{ji}^2+s_{ji}) + \frac{\tilde{\rho}_{s_{ji}}}{2}(p_{ji}^2+q_{ji}^2+s_{ji})^2\\
    &\text{subject to} \;\; \underline{v}_i \le v_i \le \overline{v}_i,\;
    \underline{v}_j \le v_j \le \overline{v}_j\\
    &\qquad\qquad \;\; -2\pi \le \theta_{i(ij)} \le 2\pi,\;
    -2\pi \le \theta_{j(ij)} \le 2\pi\\
    &\qquad\qquad \;\; -\bar{r}_{ij} \le s_{ij} \le 0,\;
    -\bar{r}_{ji} \le s_{ji} \le 0\\
    &\text{where}\;\; p_{ij},q_{ij},p_{ji},q_{ji} \; \text{are replaced by}\; \eqref{eq:pij}\text{--}\eqref{eq:qji}.
    \end{aligned}
    \label{eq:branch}
\end{equation}

Formulations for generators and buses are equivalent to~\cite{Mhanna19} except that we have an artificial variable $z$ as in~\eqref{eq:two-level} for convergence guarantees.
Since $z$ is fixed when we solve those components subproblems, generators and buses have a closed-form solution as in~\cite{Mhanna19}.

\subsection{Convergence of our ADMM iterations}
\label{subsec:convergence}

Let $N$ denote the total number of components (the total number of generators, branches, and buses) of a given power grid.
Our formulation to apply the two-level ADMM algorithm encapsulated in~\eqref{eq:acopf-admm}, where $x_i$ denotes generator and branch component variables,  $\bar{x}$ is a bus component variable,  and $z$ is an artificial variable.
For a branch component $i$, we note that its feasible region $X_i$ is nonconvex.
For generators and buses, their feasible regions are convex, as described in Section~\ref{subsec:gpu-implementation}.
Without loss of generality, we can assume that $X_i$ and $\bar{X}$ are compact.
The reason is that in practice there is a limit on power flow for each branch.
\begin{equation}
    \begin{aligned}
    &\underset{x_i \in X_i,\bar{x}\in\bar{X},z}{\text{minimize}} && \sum_{i=1}^N f_i(x_i)\\
    &\text{subject to} && Ax + B\bar{x} + z = 0, \quad
    z=0\\
    &\text{where} && x=(x_1,\dots,x_N)
    \end{aligned}
    \label{eq:acopf-admm}
\end{equation}
By taking an augmented Lagrangian of~\eqref{eq:acopf-admm} on the $z=0$ constraint and applying similar techniques to prove Theorem 1 of~\cite{SunSun21}, one can easily verify that global convergence to a stationary point of the original ACOPF problem follows.

\section{Implementation of the ADMM Algorithm on GPUs}
\label{sec:implementation}

\begin{algorithm}[t]
\begin{algorithmic}[1]
\REQUIRE{Our ADMM formulation described in Section~\ref{subsec:our-admm}}
\ENSURE{a solution $(x^*,\bar{x}^*,z^*,y^*)$}
\REPEAT
\WHILE{inner iteration not converged}
\STATE{Solve for $x$ (solutions for generators and branches).}
\STATE{Solve for $\bar{x}$ (solutions for buses).}
\STATE{Solve for $z$.}
\STATE{Update multiplier $y$.}
\ENDWHILE
\STATE{Update multiplier $\lambda$ and penalty $\beta$.}
\UNTIL{$\|z\| \nleq \epsilon$}
\end{algorithmic}
\caption{Two-level algorithm for component-based decomposition}
\label{fig:algorithm}
\end{algorithm}

Algorithm~\ref{fig:algorithm} describes the overall procedure of our ADMM algorithm to solve an ACOPF problem.
Basically, our ADMM algorithm is designed to  utilize massively parallel computing while retaining convergence guarantees: each of the solve and multiplier update routines described in lines 3--8 of Algorithm~\ref{fig:algorithm} can be executed in parallel.
We have implemented these routines entirely on GPUs, without any data transfer between CPUs and GPUs.

In this section we present our design considerations to efficiently implement an ADMM algorithm on GPUs, and we discuss the reason why we may want to use the component-based decomposition instead of the network decomposition scheme of~\cite{SunSun21}.
We then introduce the details of our implementation on GPUs for each solve and multiplier update routine.

\subsection{Design principles on efficient GPU implementation}
\label{subsec:design-principles}


In~\cite{Hogg16} the authors presented a GPU-based linear system solver for symmetric indefinite system of equations. Its application in Ipopt to ACOPF problems~\cite{Tasseff19}, however, has shown a much slower computation time than the times obtained from using existing CPU-based linear solvers.
In~\cite{Swirydowicz21} the authors compared performance between a variety of existing CPU- and GPU-based linear solvers (including Nvidia's cuSOLVER) on linear systems of equations that originate from ACOPF problems.
Their findings demonstrate that no significant acceleration has been observed by using GPU-based linear solvers.

In this regard, if we were to use the network decomposition of~\cite{SunSun21}, we would have to solve a number of ACOPF subproblems on GPUs.
As long as we use the existing optimization algorithms, such as Ipopt, computation time  will significantly depend on that of linear system solves, which will in turn show degraded performance on GPUs, as observed in~\cite{Tasseff19,Swirydowicz21}.

These observations in the literature lead us to look for other approaches that do not heavily rely on accelerating the solve time of a sparse linear system of equations.
Component-based decomposition provides favorable problem structures that can be massively exploited via the GPU's parallel computing capability.
Generators and buses subproblems in lines 3--4 of Algorithm~\ref{fig:algorithm} have a closed-form solution as well as the solves and updates in lines 5--8, as described in Section~\ref{subsec:gpu-implementation}.
Therefore, their solutions can be directly and efficiently computed in parallel on GPUs by employing a large number of threads.
The only exceptions are branch subproblems, which are nonlinear nonconvex optimization problems as described in~\eqref{eq:branch}.
In this case we employ our recently developed GPU-based batch nonlinear programming solver~\cite{ExaTron}.

From the standpoint of parallel computing, another benefit of component-based decomposition is that it is highly scalable.
The size of each subproblem is independent of the size of the grid network of a given ACOPF problem.
Only the number of subproblems increases in this case.
This independence implies that we can easily scale up the ADMM algorithm by employing more computational resources, such as the number of GPUs, to deal with  additional subproblems in parallel.

\subsection{Our GPU-based implementation}
\label{subsec:gpu-implementation}

Mathematically, the subproblems for generators (line 3 of Algorithm~\ref{fig:algorithm}) and for $z$ (line 5 of Algorithm~\ref{fig:algorithm}) take the following form:
\begin{equation}
    \begin{aligned}
    &\underset{l \le x \le u}{\text{minimize}} \; \frac{1}{2}x^TQx - c^Tx\\
    (\Rightarrow)\; & x^*=\max(l, \min(u, Q^{-1}c)),
    \end{aligned}
    \label{eq:qp-bound}
\end{equation}
where $Q$ is strongly convex.\footnote{We abuse the notation $x$ here. In~\eqref{eq:qp-bound}, $x$ is not the same as $x$ in line 3 of Algorithm~\ref{fig:algorithm}.}
In the case of generators, entries for $Q$ come from their quadratic objective function and penalty terms.
For the $z$ variable, $Q$ is formed from their penalty terms.
We note that $(l,u)=(-\infty,\infty)$ when we solve for the $z$ variable.

The subproblems for buses have the following form where $Q$ (of their penalty terms) is strongly convex as well.
\begin{equation}
    \begin{aligned}
    &\underset{x}{\text{minimize}} \; \frac{1}{2}x^TQx - c^Tx\\
    &\text{subject to} \; Ax=b \qquad\qquad (\mu)\\
    (\Rightarrow)\; & \mu^* = (AQ^{-1}A^T)^{-1}(AQ^{-1}c-b)\\
    & x^*= Q^{-1}(c-A^T\mu^*)
    \end{aligned}
    \label{eq:qp-bus}
\end{equation}

In~\eqref{eq:qp-bus}, $\mu$ denotes a multiplier for constraint $Ax=b$.
The matrix $A$ has two rows corresponding to the real and reactive power flow equations~\eqref{eq:pg}--\eqref{eq:qg}, respectively.
They are linearly independent of each other; therefore, $AQ^{-1}A^T$ is nonsingular, and $\mu^*$ is well defined.

For the multiplier updates, the following update rules are used:
\begin{equation}
    \begin{aligned}
    y &\leftarrow y + \rho(Ax + B\bar{x} + z)\\
    \lambda &\leftarrow \Pi_{[\underline{\lambda},\overline{\lambda}]}\left(\lambda + \beta z\right),
    \end{aligned}
\end{equation}
where $\Pi_{[\underline{\lambda},\overline{\lambda}]}$ is a projection operator onto their lower and upper bounds, $\underline{\lambda}$ and $\overline{\lambda}$, respectively.

Since the aforementioned subproblems and multiplier update rules have a closed-form solution, their implementation on GPUs is straightforward:
we launch the same number of threads on GPUs as the number of elements in $x, y$, and $\lambda$.
Each thread updates its corresponding element in parallel.

For branch subproblems~\eqref{eq:branch}, we employ our GPU-based batch nonlinear programming solver~\cite{ExaTron}. 
Each nonconvex nonlinear optimization branch subproblem is solved via a trust-region Newton's algorithm based on a preconditioned conjugate gradient method~\cite{LinMore99}.
Nonconvexity of the branch problem is handled via detecting and following a negative curvature direction during the conjugate gradient iteration~\cite{Steihaug83}.
A batch of such subproblems are solved on GPU in parallel by launching the same number of thread blocks (consisting of a group of threads) as the number of branch subproblems.
Each thread block solves a specific branch subproblem.

A distinguishing feature of our GPU-based batch solver is that it operates entirely on GPUs without requiring data transfers between the host (CPU) and the device (GPU) during its operation.
Since the size of each branch problem is very small (involving 6 variables only), its computational performance heavily depends on the time for memory access.
If our solver were to require data transfer between the host and the device, it would significantly slow the computation time.
Our solver avoids such expensive memory transfer by implementing the entire optimization algorithm on GPUs.


\section{Numerical Results}
\label{sec:exp}

In this section we demonstrate the computational performance of our GPU-based ADMM solver described in Sections~\ref{sec:formulation}--\ref{sec:implementation} over large power grids with up to a 70,000 bus system.
Our solver has the ability to solve ACOPF from both cold start and warm start.
In Section~\ref{subsec:single-period} we  present experimental results of solving ACOPF from cold  start , where we compare its performance with Ipopt~\cite{Wachter06}.
We then demonstrate the performance of our warm-start capability for rapidly recomputing optimal solutions upon load changes as we move forward in time.


\subsection{Experiments setting}
\label{subsec:setting}

\begin{table}[t]
\centering
\caption{Data and parameters for experiments}
\label{tbl:data}
\begin{tabular}{|r|r|r|r|r|r|}
    \hline
    \multicolumn{1}{|c|}{Data} 
    & \multicolumn{1}{c|}{\# Generators}
    & \multicolumn{1}{c|}{\# Branches}
    & \multicolumn{1}{c|}{\# Buses}
    & \multicolumn{1}{c|}{$\rho_{pq}$}
    & \multicolumn{1}{c|}{$\rho_{va}$}\\\hline
    1354pegase & 260 & 1,991 & 1,354 & 1e1 & 1e3\\
    2869pegase & 510 & 4,582 & 2,869 & 1e1 & 1e3\\
    9241pegase & 1,445 & 16,049 & 9,241 & 5e1 & 5e3\\
    13659pegase & 4,092 & 20,467 & 13,659 & 5e1 & 5e3\\
    ACTIVSg25k & 4,834 & 32,230 & 25,000 & 3e3 & 3e4\\
    ACTIVSg70k & 10,390 & 88,207 & 70,000 & 3e4 & 3e5\\
    \hline
\end{tabular}
\end{table}

The implementation of our GPU-based ADMM solver has been written in Julia@v1.6.3 using CUDA.jl@v3.4.2.
To compare the performance between our solver and Ipopt@v3.13.4, we used PowerModels.jl@v0.18.3~\cite{Carleton18} to solve ACOPF on CPUs using Ipopt with MA57 being set as its linear solver.
In this case we directly solve the ACOPF formulation~\eqref{eq:acopf} using Ipopt.
Since PowerModels.jl automatically tightens angle difference constraints when it thinks they are too loose, we disabled them from PowerModels.jl for a fair comparison.\footnote{A function call to constraint\_voltage\_angle\_difference() has been commented out in the build\_opf() function.}
Table~\ref{tbl:data} shows a list of power grid files retrieved from MATPOWER~\cite{Zimmerman11} for our experiments.
We performed experiments on Nvidia's Quadro GV100 for our GPU solver and on Intel Xeon 6140 CPU@2.30GHz for Ipopt.

Our ADMM algorithm involves a number of parameters such as multipliers and penalty terms on consensus constraints.
All experiments in this section were performed by using zero as initial values for multipliers and fixed values for penalty terms, as shown in Table~\ref{tbl:data}.
The penalty term $\rho_{pq}$ is for real and reactive power generation and power flow, and $\rho_{va}$ is for voltage magnitude and angle.
Their values are fixed throughout the ADMM iterations.
We note that large penalty values lead to less weight on the objective value, potentially causing a larger optimality gap.
To reduce this adverse effect, we scaled the objective value for the 70k case by multiplying it by 2.
The termination conditions for our ADMM algorithm are the same as~\cite{SunSun21}, where we set the maximum number of outer and inner iterations to 20 and 1,000, respectively.
In all our experiments the termination conditions have been satisfied before reaching the maximum number of outer iterations.

As a solution for our ADMM algorithm, we use real and reactive power generation from generator problems and voltage values from bus problems.
Instead of using power flow values from branch problems, they are recomputed from the voltage values of bus problems for consistency.
In this case, consensus errors between bus and branch problems could propagate to line limit constraints, resulting in higher line limit violations than computed by~\eqref{eq:branch}.
To reduce this type of error propagation, we slightly tighten the line limit by allowing it to use up to 99\% of its capacity.


\subsection{Performance of solving ACOPF from cold start}
\label{subsec:single-period}

We measured the computational performance of cold start of our solver and Ipopt by solving ACOPF over power grids described in Table~\ref{tbl:data}.
For both solvers we set the initial values for real and reactive power generation and voltage magnitude to the medium of their lower and upper bounds, respectively.
Voltage angles were initialized to zero, and the reference voltage angle was fixed to zero.
The initial values for power flows were computed based on the initial voltage values.

Table~\ref{tbl:single-period} presents the computational performance of our solver and Ipopt.
The ADMM Iterations column shows the cumulative number of inner-level ADMM iterations of our solver.
$\|c(x)\|_\infty$ and $\frac{|f-f^*|}{f^*}$ represent the maximum constraint violation and the relative gap of the objective value at a solution obtained from our solver, respectively.
The relative gap was measured by using the objective value $f^*$ from Ipopt.

The experimental results in Table~\ref{tbl:single-period} demonstrate that our GPU-accelerated solver was able to quickly find a solution of good quality from cold start.
For pegase cases, the maximum constraint violations were on the order of $10^{-3}$ or $10^{-4}$, while the relative objective gap was less than 0.1\%.
For 25k and 70k cases, the constraint violations deteriorated but were still on the order of $10^{-2}$ with the relative objective gap being less than $2.5\%$.
We note that the termination conditions of ADMM algorithms typically become less strict as the number of consensus constraints increases.
In our case the number of such constraints is on the order of $10^5$ for the two largest cases.

To the best of our knowledge, our solver is the first GPU-based solver that shows such competitive computational performance as indicated in Table~\ref{tbl:single-period}, compared with Ipopt.

\begin{table}[t]
\centering
\caption{Performance of solving ACOPF from cold-start}
\label{tbl:single-period}
\begin{tabular}{|r|r|r|r|r|r|}
    \hline
    \multicolumn{1}{|c|}{\multirow{2}{*}{Data}} 
    & \multicolumn{1}{c|}{ADMM}
    & \multicolumn{2}{c|}{Time (secs)}
    & \multicolumn{2}{c|}{Solution Quality}\\\cline{3-6}
    & \multicolumn{1}{c|}{Iterations} 
    & \multicolumn{1}{c|}{ADMM} & \multicolumn{1}{c|}{Ipopt} 
    & \multicolumn{1}{c|}{$\|c(x)\|_\infty$}
    & \multicolumn{1}{c|}{$\frac{|f-f^*|}{f^*}$}\\\hline
    1354pegase   &   823 &   1.99 &   2.44 & 1.23e-03 & 0.05\%\\
    2869pegase   & 1,230 &   4.19 &   6.09 & 3.64e-04 & 0.03\%\\
    9241pegase   & 1,372 &   7.95 &  50.80 & 1.12e-03 & 0.08\%\\
    13659pegase  & 1,529 &   8.70 & 131.12 & 1.25e-03 & 0.05\%\\
    ACTIVSg25k   & 3,307 &  36.05 & 118.64 & 1.21e-02 & 0.09\%\\
    ACTIVSg70k   & 2,897 &  69.81 & 469.03 & 1.52e-02 & 2.20\%\\
    \hline
\end{tabular}
\end{table}


\subsection{Performance of solving ACOPF from warm start}
\label{subsec:warm-start}


To measure the computational performance of warm start, we performed experiments  solving ACOPF over a time horizon, where each time period is solved by warm starting from a solution obtained from the previous time period.
The time horizon has 30 time periods with each period representing one minute.
Its load profile was generated by interpolating an hourly real-time system demand data from ISO New England~\cite{ISOdata} into minutes.
During the 30-minute time horizon, the load changes up to 5\% from its starting value.
We note that when we warm start from a given point, we take account of the ramp rates of generators so that it holds $|p_{g_i,t+1}-p_{g_i,t}| \le r_g$ for $t=1,\dots,29$, where $r_g$ is set to 2\% of generator $g$'s upper limit of real power generation.
To warm start Ipopt, we set \texttt{warm\_start\_init\_point} option to \texttt{yes}.

\begin{figure}[h]
\centering
\includegraphics[scale=.65]{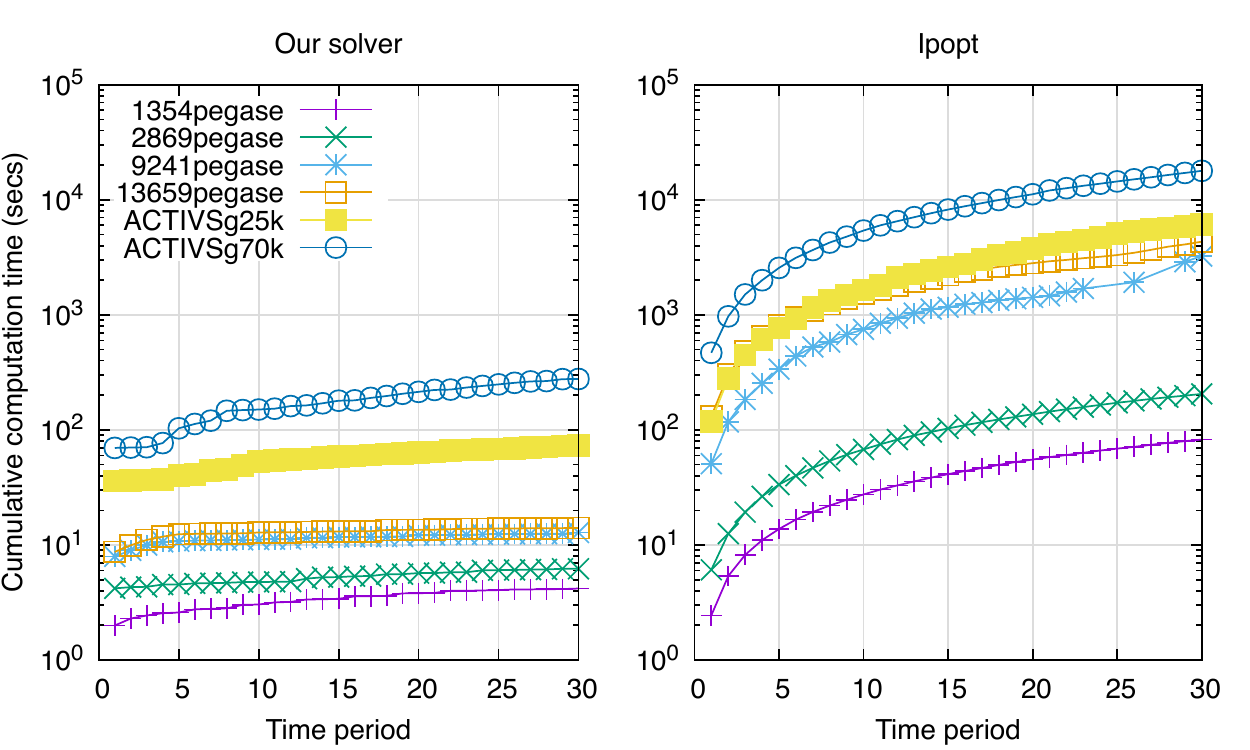}
\caption{Cumulative computation time of warm start}
\label{fig:warmstart-time}
\end{figure}

Figure~\ref{fig:warmstart-time} presents the cumulative computation time of our solver and Ipopt over the time horizon, where we solved the first time period with cold start and we warm started the subsequent periods. 
As the figure illustrates, our warm start shows a significantly accelerated computational performance: for pegase data in most periods it took less than 1 second to solve; for the 25k case the computation time was usually in  1--2 seconds, with some outliers taking about 4 seconds; for 70k mostly it took less than 10 seconds for each period, with some exceptions  taking 25 seconds.
We highlight that our GPU solver solves the entire time horizon of ACTIVSg70k faster than does Ipopt only for the first single time period.

Compared with our solver, computational gains were not observed with Ipopt when we warm started it.
The cumulative computation time linearly increased as we moved along the time horizon, implying that Ipopt was not taking advantage of warm start.
We initially tried with additional options such as preventing the given point from being pushed into the interior as shown in~\cite{Kim20}; however, we found that those options did not provide any computational benefits in our case.
We think that the reason is mainly that our load fluctuates more over time than in~\cite{Kim20}, where the load changes were less than 2\% over 30 minutes.
In addition, we note that, as shown in Figure~\ref{fig:warmstart-time}, Ipopt failed to solve the 9241pegase case for some time periods because of numerical issues.

\begin{figure}[t]
\centering
\includegraphics[scale=.65]{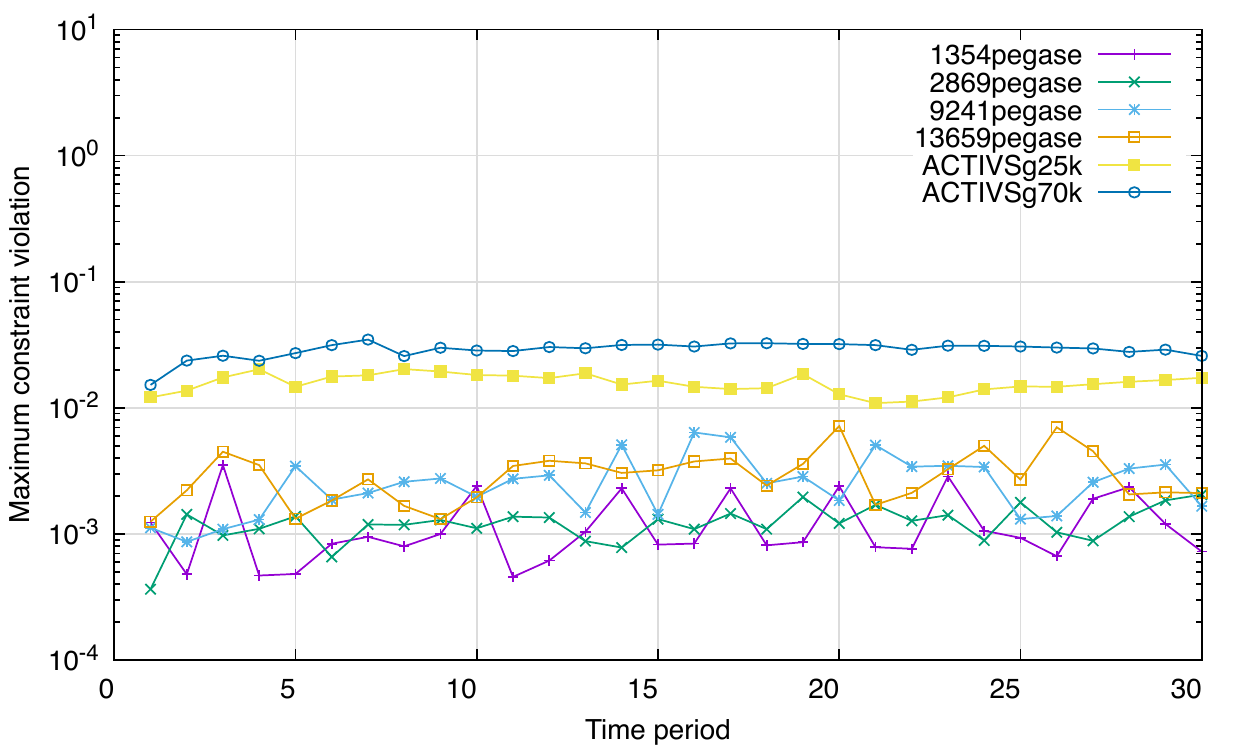}
\caption{Maximum constraint violation of warm start}
\label{fig:warmstart-viol}
\end{figure}

The maximum constraint violations along the time horizon are presented in Figure~\ref{fig:warmstart-viol}.
The violations stayed on the order of values similar to those we obtained with cold start as listed in Table~\ref{tbl:single-period}.
No significant deterioration was observed as we moved forward in time.
Since we solve the original ACOPF~\eqref{eq:acopf} using ADMM until convergence (unlike solving some simplified formulation), the solution quality is expected to be maintained as with cold-start.

\begin{figure}[t]
\centering
\includegraphics[scale=.65]{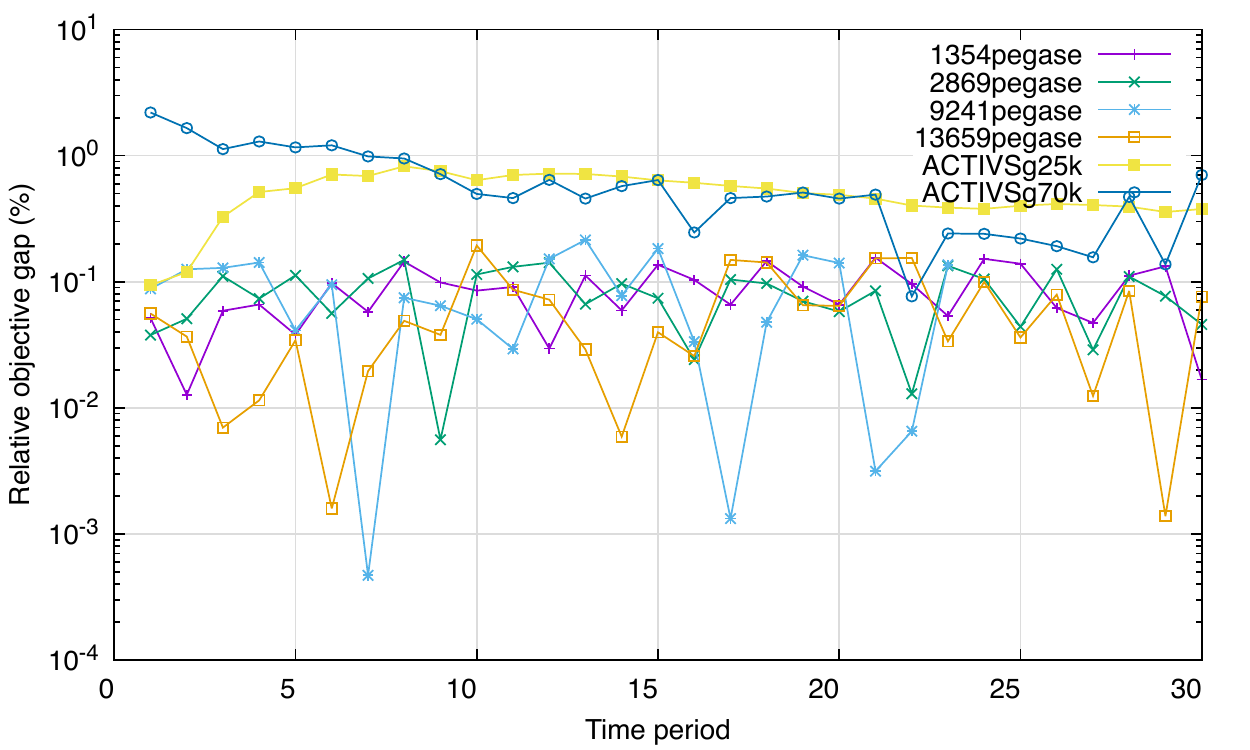}
\caption{Relative objective gap of warm start}
\label{fig:warmstart-obj}
\end{figure}

Figure~\ref{fig:warmstart-obj} depicts the relative objective gap for each time period.
Similar to the results of Figure~\ref{fig:warmstart-viol}, the quality of the relative objective gap was maintained at a similar level to that obtained with cold start.
We note that after the seventh period all of the relative objective gaps became less than 1\%.

Overall, our warm start demonstrated a substantial improvement of computation time compared with cold start while maintaining the solution quality.
The computation time for pegase data was almost in real time, taking less than 1 second for each time period.
For larger data, however, it took more time (i.e., 30 seconds per time period).

\section{Conclusion}
\label{sec:conclusion}

We have presented the first fully GPU-based scalable optimization method based on ADMM for solving an ACOPF with convergence guarantees.
To rapidly solve a batch of small nonlinear nonconvex subproblems, we developed augmented Lagrangian techniques that enable us to reformulate those subproblems into a batch of bound-constrained optimization problems, in order to exploit our GPU-based batch solver ExaTron~\cite{ExaTron}.
The numerical experiments report a significant reduction in solution wall-clock times with our GPU ADMM solver, as compared with the Ipopt solution times, over large grids having up to a 70,000 bus system.
Moreover, the warm-start capability of our solver demonstrated a significantly accelerated computational performance: less than a second to solve (a real-time performance) for up to a 13k bus system and less than 30 seconds for a 70k system, making it particularly suitable for rapidly tracking optimal set points of generators under frequent load and generation fluctuations.

Our ADMM solver has considerable potential for improvement.
The distributed and scalable nature of our solution method allows us to further accelerate the computation time by utilizing multiple GPUs.
As observed in~\cite{Mhanna19,SunSun21}, penalty terms of the ADMM algorithm could significantly affect its computation time until convergence.
An automatic penalty selection scheme based on machine learning techniques, such as reinforcement learning, may greatly reduce the number of ADMM iterations until convergence.



%
\bibliographystyle{IEEEtran}
\bibliography{pscc2022}

\vspace{3mm}
{\footnotesize
\noindent\fbox{\parbox{0.47\textwidth}{
The submitted manuscript has been created by UChicago Argonne, LLC, Operator of Argonne National Laboratory (``Argonne''). Argonne, a U.S. Department of Energy Office of Science laboratory, is operated under Contract No. DE-AC02-06CH11357. The U.S. Government retains for itself, and others acting on its behalf, a paid-up nonexclusive, irrevocable worldwide license in said article to reproduce, prepare derivative works, distribute copies to the public, and perform publicly and display publicly, by or on behalf of the Government. The Department of Energy will provide public access to these results of federally sponsored research in accordance with the DOE Public Access Plan (http://energy.gov/downloads/doe-public-access-plan).}
}
}

\end{document}